\documentclass[preprint,12pt]{elsarticle}

\usepackage{graphicx}%
\usepackage{multirow}%
\usepackage{amsmath,amssymb,amsfonts}%
\usepackage{amsthm}%
\usepackage{mathrsfs}%
\usepackage[title]{appendix}%
\usepackage{xcolor}%
\usepackage{textcomp}%
\usepackage{manyfoot}%
\usepackage{booktabs}%
\usepackage{algorithm}%
\usepackage{algorithmicx}%
\usepackage{algpseudocode}%
\usepackage{listings}%
\usepackage{cleveref}
\usepackage{tikz}

\newtheorem{theorem}{Theorem}[section]
\newtheorem{definition}[theorem]{Definition}
\newtheorem{lemma}[theorem]{Lemma}
 
\newtheorem{corollary}[theorem]{Corollary}
\newtheorem{example}[theorem]{Example}
\newtheorem{proposition}[theorem]{Proposition}

\begin{document}

\begin{frontmatter}

\title{On the structure and spectra of an induced subgraph of essential ideal graph of $\mathbb{Z}_{n}$}





\author[1]{Jamsheena P}
\affiliation[1]{organization={Department of Mathematics, National Institute of Technology Calicut},
            addressline={NIT Calicut Campus(PO)},
            city={Kozhikode},
            postcode={673601},
            state={Kerala},
            country={India}}


 \author[1]{ Chithra A V}

\begin{abstract}
Let $R$ be a commutative ring with unity. The essential ideal graph   $\mathcal{E}_R$  of $R$ is a graph in which the vertex set comprises of set of all nonzero proper ideals of $R$ and two vertices $I$ and $K$ are adjacent if and only if $I+K$ is an essential ideal. In this paper, we discuss the structure of an induced subgraph of the essential ideal graph of the ring $\mathbb{Z}_{n}$ as a $\mathscr{G}$-generalized join graph and thereby completely determine the structure of $\mathcal{E}_{\mathbb{Z}_{n}}$. Also, we prove a characterization of $\mathcal{E}_{\mathbb{Z}_{n}}$ to be Laplacian integral in terms of the vertex-weighted Laplacian matrix of annihilating ideal graph of $\mathbb{Z}_{n}$  for $n= \prod_{i=1}^k p_i$. Further, we discuss the eigenvalues of various matrices like adjacency matrix, Laplacian matrix, signless Laplacian matrix, and normalized Laplacian matrix of the induced subgraph of the essential ideal graph of $\mathbb{Z}_{n}$. Finally, we obtain the upper bounds of spectral radius and algebraic connectivity of  $\mathcal{E}_{\mathbb{Z}_{n}}$ and compute the values of $n$ for which these bounds are attained.
\end{abstract}



\end{frontmatter}
\section{Introduction}
Let $\Gamma=(V, E)$ be an undirected finite simple graph having vertex set $V(\Gamma)=\{ u_1,u_2, \cdots, u_n\}$ and edge set $E$ consisting of unordered pairs of vertices. If the vertex $u_i$ is adjacent to a vertex $u_j$ for $j\ne i$ in $\Gamma$, then we write $u_i\sim u_j$. Also, $N_\Gamma(u)=\{v\in V(\Gamma): v\sim u\ \text{in}\ \Gamma\}$ and $deg(u)= |N_\Gamma(u)|$ denote
the set of neighbors and degree of the vertex $u$, respectively. 
 For a graph $\Gamma$, a partition of vertices   $V(\Gamma)= V_1\cup V_2 \cup \cdots \cup V_k$ is said to be an \textit{equitable partition} if each vertex in $V_i$ has the same number of neighbors in $V_j$ for any $i,j \in \{1, 2, \cdots, k \}$. The \textit{join} $\Gamma_1 \vee \Gamma_2$ of two graphs $(\Gamma_1, n_1)$ and $\Gamma_2, n_2)$, is a graph obtained from $\Gamma_1$ and $\Gamma_2$ by joining each vertex of $\Gamma_1$ to every vertex of $\Gamma_2$.
The \textit{adjacency matrix} is the $n\times n$ matrix $A(\Gamma)=(a_{ij})$, where $a_{ij} =\begin{cases}
 1 & \textit{if}  \ v_i\sim v_j \\ 
 0 &   \textit{if} \  v_i\nsim v_j \\
\end{cases}$ in $\Gamma$. The matrices $L(\Gamma)= D(\Gamma)-A(\Gamma)$ and $Q(\Gamma)= D(\Gamma)+ A(\Gamma)$, $D(\Gamma)$ is the diagonal matrix of vertex degrees, are called the \textit{Laplacian matrix} and \textit{signless Laplacian matrix} of $\Gamma$, respectively. Since $L(\Gamma)$ is real, symmetric, and positive semi-definite, all its eigenvalues called Laplacian eigenvalues are real and non-negative. The signless Laplacian matrix  $Q(\Gamma)$
\cite{cvetkovic2009towards}, is also a real, symmetric, and positive semi-definite matrix having non-negative eigenvalues. This was further studied in \cite{cvetkovic2010towards,cvetkovic2010towardsIII}. The \textit{normalized Laplacian matrix} of $\Gamma$, 
 is defined as $\mathscr{L}(\Gamma)= D^{-\frac{1}{2}}L(\Gamma)D^{\frac{1}{2}}$, where $D^{-\frac{1}{2}}$ represents the diagonal matrix having $i$-th diagonal entry $\frac{1}{\sqrt{deg(u_i)}}$ for $1\le i\le n$ \cite{chung1997spectral}.

Let $\sigma(C)$ denote the multiset of all eigenvalues of the square matrix $C$, called the spectrum of $C$. The \textit{Laplacian spectrum} of $\Gamma$ is the multiset of all Laplacian eigenvalues $\mu_1> \mu _2> \cdots >\mu_{l}$   with respective multiplicities   $s_1,s_2,\cdots, s_l$ 
and is denoted by $\sigma_L(\Gamma)= \begin{pmatrix}  
    \mu_1 & \mu_2 & \cdots & \mu_l \\
    s_1 & s_2 & \cdots & s_l \\
    \end{pmatrix}$. In a similar way, $\sigma_A(\Gamma)$, $\sigma_Q(\Gamma)$ and $\sigma_{\mathscr{L}}(\Gamma)$ denote the \textit{adjacency, signless Laplacian} and \textit{normalized Laplacian spectrum} of the graph $\Gamma$, respectively. A graph is said to be \textit{Laplacian integral} if all its Laplacian eigenvalues are integers. The \textit{algebraic connectivity} of a graph $\Gamma$ is the second smallest eigenvalue $\mu_{n-1}$ of $\Gamma$ and is denoted as $a(\Gamma)$. The largest eigenvalue $\mu_1$ of $L(\Gamma)$ is called the \textit{spectral radius} of the graph $\Gamma$ and is denoted as $b(\Gamma)$. It is proved that  $b(\Gamma)= n-a(\overline{\Gamma})$ \cite{fiedler1973algebraic}, where $\overline{\Gamma}$ denotes the complement graph of $\Gamma$. The \textit{vertex connectivity}, $\kappa(\Gamma)$, of a graph $\Gamma$ is the minimum number of vertices needed to be removed to disconnect the graph. For a complete graph $K_n$ on $n\ge 2$ vertices,  $\kappa(K_n)=n-1$.

Let $R$ be a commutative ring with nonzero unity. An element $z\in R$ is said to be a zero divisor of $R$ whenever there exists a nonzero element $w\in R$ such that $zw=0$. An ideal $I$ of a ring $R$ which has a nonzero intersection with every other nonzero ideal of $R$ is called an \textit{essential ideal}.  An ideal $I$ of a ring $R$ is said to be an \textit{annihilating ideal} if a nonzero ideal $K$ of $R$ exists such as $IK=0$.

The study of \textit{zero divisor graphs} introduced by Beck \cite{beck1988coloring} was the origin of the era of research on graphs from algebraic structures. It was further modified by Anderson and Livingston \cite{ANDERSON1999434}. Following this, a large variety of graphs were defined and studied in the literature. But, in ring theory, the structure of a ring is exactly linked to how ideals behave rather than the individual elements themselves. Based on this, research was continued on graphs from rings by taking the ideals of a ring as vertices.
 One can refer \cite{alilou2016sum, behboodi2011annihilating, ye2012co} for details. In 2011 \cite{behboodi2011annihilating, behboodi2011annihilatingII}, Behboodi and Rakeei defined and studied the annihilating ideal graph, $\mathbb{AIG}(R)$, of a commutative ring $R$. It is an undirected graph with a vertex set consisting of a set of all nonzero annihilating ideals and two distinct ideals $I$ and $K$ are adjacent if and only if $IK=0$.  This was further studied in \cite{aalipour2014classification,chelvam2015connectivity, curtis2018classifying}.  In 2012, Ye and Wu \cite{ye2012co} defined \textit{comaximal ideal graph} of a commutative ring with vertices as the proper ideals which are not contained in $J(R)$, the Jacobson radical of $R$ and two distinct vertices $I$ and $K$ are adjacent whenever $I+K=R$. In $2016$, Azadi et al.\cite{azadi2016comaximal} investigated the planarity and perfectness of the comaximal ideal graph.

Recently, in 2018 Amjadi \cite{amjadi2018essential} introduced and studied a large class of graphs named \textit{essential ideal graph}  $\mathcal{E}_R$   which contains comaximal ideal graph as a subgraph, for a commutative ring $R$. He proved that $\mathcal{E}_R$ is a connected graph with diameter $\le 3$. In this paper, we continue the study of essential ideal graphs for the integers modulo ring, $\mathbb{Z}_{n}$. We identify that the graphical structure of the induced subgraph of nonessential ideals of the ring $\mathbb{Z}_{n}$ is the $\mathscr{G}$-generalized join of certain null graphs. Also, it is proved that the graph  $\mathscr{G}$ of the $\mathscr{G}$-generalized join
can be identified as the annihilating ideal graph of the ring  $\mathbb{Z}_{n}$, where $n$ is the product of distinct primes.

This paper has been arranged subsequently. In Section $2$, we give some definitions and theorems which are used in this paper. In Section $3$, we determine the graphical structure of the induced subgraph $\mathcal{E}_{\mathbb Z_n}(\mathscr{U})$, $\mathscr{U}$ being the set of nonessential ideals of $\mathbb Z_n$ and thereby determine the structure of $\mathcal{E}_{\mathbb Z_n}$ completely. In Section $4$, we compute the spectrum of various matrices like adjacency, Laplacian, signless Laplacian, and normalized Laplacian of the induced subgraph $\mathcal{E}_{\mathbb Z_n}(\mathscr{U})$. Also, we obtain a characterization of $\mathcal{E}_{\mathbb Z_n}$ to be Laplacian integral by the vertex-weighted Laplacian matrix. In Section $5$, we discuss the spectral properties such as spectral radius and algebraic connectivity of the essential ideal graph $\mathcal{E}_{\mathbb Z_n}$.  We calculate the upper bounds for the spectral radius and algebraic connectivity of $\mathcal{E}_{\mathbb Z_n}$.


\section{Preliminaries}
In this section, we list some definitions and preliminary theorems needed in subsequent sections.
\begin{lemma}\cite{zhang2006schur} \label{detMATX}
	Let $S= \begin{pmatrix}
	 M & N \\
	 P & Q  \\
	\end{pmatrix}$, where $M,N,P,Q$ be matrices and $Q$ be non-singular. Then $\det S =  \det Q \times \det (M- NQ^{-1} P)$. 
 \end{lemma}
\begin{theorem}  \label{chr. essentl.idl Zn} \cite{jamsheena2023adjacency}
   Let $n= p_{1}^{m_1}p_{2}^{m_2} \cdots p_{k}^{m_k}$ where $p_{1}<p_{2}< \cdots < p_{k}$ are primes, and $m_j$ is a non-negative integer for $1 \le j \le k$. Any nonzero ideal $I= \langle p_{1}^{r_1}p_{2}^{r_2} \cdots p_{k}^{r_k} \rangle $ of $\mathbb {Z}_n$ is essential if and only if $r_j \ne m_{j}$ for any $j$. 
\end{theorem}
\begin{theorem} \label{Strctr thm Zn} \cite{jamsheena2021structural}
For any composite integer $n=p_{1}^{\alpha_{1}}p_{2}^{\alpha_{2}} \cdots p_{k}^{\alpha_{k}}>1$, $(k,\ \alpha_{j})\in \mathbb{N} $,  $(k,\alpha_{1})\neq (1,1)$, $p_{j}$'s are distinct primes ($1\le j\le k$), \\
 the essential ideal graph $\mathcal{E}_{\mathbb{Z}_{n}}\cong K_{m} \vee H $, where $K_{m}$ is a complete graph of order $m= (\displaystyle \prod_{j=1}^{k}\alpha_{j}-1)$ and $H$ is a $k$-partite graph .
\end{theorem}
\begin{theorem}\label{LcharPlnml join}\cite{mohar1991graph}
    Let $(\Gamma_1, n_1)$ and $(\Gamma_2,n_2)$ be two vertex disjoint graphs. Then the Laplacian characteristic polynomial of $\Gamma_1\vee \Gamma_2$ is given by 
    \begin{equation*}
        P_L(\Gamma_1\vee \Gamma_2, x)= \frac{x(x-n_1-n_2)}{(x-n_1)(x-n_2)} P_L(\Gamma_1,x-n_2) P_L(\Gamma_2, x-n_1).
    \end{equation*}
\end{theorem}
\begin{definition}\cite{schwenk2006computing}
    Let $\mathscr{G}$ be a graph on $k$ vertices $v_1,v_2,\cdots,v_k$ and $\Gamma_1, \Gamma_2,\\
    \cdots,\Gamma_k$ be pairwise disjoint graphs. Then, the $\mathscr{G}$- generalized join graph $\mathscr{G}[\Gamma_1, \Gamma_2,\cdots,\Gamma_k]$ of $\Gamma_1, \Gamma_2,\cdots,\Gamma_k$ is the graph formed by replacing each vertex $v_i$ of $\mathscr{G}$ by the graph $\Gamma_i$ and joining each vertex of $\Gamma_i$ to every vertex of $\Gamma_j$ whenever $v_i\sim v_j$ in $\mathscr{G}$.
\end{definition}

\begin{theorem}\label{spec G-join} \cite{cardoso2013spectra, wu2014signless} 
    Let $\mathscr{G}$ be a graph on $k$ vertices $v_1,v_2,\cdots,v_k$ and let $\Gamma_i$ be pairwise disjoint $r_i$-regular graphs with order $n_i$, for $1\le i\le k$. Then 
\begin{itemize}
   \item \label{adj G-join}The adjacency spectrum of $\mathscr{G}[\Gamma_1, \Gamma_2,\cdots,\Gamma_k]$ is given by\\
    \begin{equation*}
        \sigma_A(\mathscr{G}[\Gamma_1, \Gamma_2,\cdots,\Gamma_k]= \bigg(\bigcup_{i=1}^{k}(\sigma_A(\Gamma_i)\backslash \{r_i\})\bigg)\cup \sigma(C_A(\mathscr{G})),
    \end{equation*}where $C_A(\mathscr{G})$ is a symmetric matrix obtained as \begin{equation*}
       C_A(\mathscr{G})=(c_{ij})= \begin{cases}
            r_i & \text{if} \  i=j, \\
            \sqrt{n_i n_j}& \text{if} \ v_i \sim v_j \ \text{in} \ \mathscr{G}, \\
            0 & \text{otherwise},
        \end{cases}
    \end{equation*} and $\sigma_A(\Gamma_i)\backslash\{r_i\}$ means that one copy of the eigenvalue $r_i$ is removed from multiset $\sigma_A(\Gamma_i)$.
 \item  \label{Lap  G-join} The Laplacian spectrum of $\mathscr{G}[\Gamma_1, \Gamma_2,\cdots,\Gamma_k]$ is given by
      \begin{equation*}
        \sigma_L(\mathscr{G}[\Gamma_1, \Gamma_2,\cdots,\Gamma_k]= \bigg(\bigcup_{i=1}^{k}(N_i+(\sigma_L(\Gamma_i)\backslash \{0\}))\bigg)\bigcup \sigma(C_L(\mathscr{G})),
    \end{equation*}where
    $N_j=\begin{cases}
        \sum_{v_i\sim v_j}n_i & \text{if}\ N_\mathscr{G}(v_j)\ne\phi \\
        0 & \text{otherwise} \\
    \end{cases}$ and $C_L(\mathscr{G})$ is a symmetric matrix obtained as \begin{equation*}
    C_L(\mathscr{G})=(c_{ij})=\begin{cases}
           N_i &  \text{if} \  i=j, \\
           - \sqrt{n_i n_j}  & \text{if}\ v_i\sim v_j \ \text{in} \ \mathscr{G}, \\
            0 & \text{otherwise}.\\
       \end{cases} 
    \end{equation*}
\item The signless Laplacian spectrum of  $\mathscr{G} [\Gamma_1, \Gamma_2,\cdots,\Gamma_k]$ is given by
      \begin{equation*}
        \sigma_{Q}(\mathscr{G}[\Gamma_1, \Gamma_2,\cdots,\Gamma_k]= \bigg(\bigcup_{i=1}^{k}(N_i+(\sigma_Q(\Gamma_i)\backslash \{2r_i\}))\bigg)\bigcup \sigma(C_Q(\mathscr{G})),
    \end{equation*}where
\begin{equation*}
        C_Q(\mathscr{G})=(c_{ij})= \begin{cases}
            2r_i+N_i, & \text{if} \  i=j,\\
            \sqrt{n_i n_j}, & \text{if} \ v_i \sim v_j \ \text{in} \ \mathscr{G}, \\
            0, & \text{otherwise}.
\end{cases}
\end{equation*} 
\item  The normalized Laplacian spectrum of  $\mathscr{G} [\Gamma_1, \Gamma_2,\cdots,\Gamma_k]$ is given by  
\begin{equation*}
\begin{split}
 \sigma_{\mathscr{L}}(\mathscr{G}[\Gamma_1, \Gamma_2,\cdots,\Gamma_k]= & \bigg(\bigcup_{i=1}^{k}\bigg(\frac{N_i}{r_i+N_i}+\frac{r_i}{r_i+N_i}(\sigma_{\mathscr{L}}(\Gamma_i)\backslash \{0\})\bigg)\bigg)\bigcup 
 \\
 & \sigma( C_{\mathscr{L}}(\mathscr{G})),
\end{split}
\end{equation*}where  
\begin{equation*}
       C_{\mathscr{L}}(\mathscr{G})=(c_{ij})= \begin{cases}
            \frac{N_i}{r_i+N_i}, & \text{if} \  i=j, \\
            \sqrt{\frac{n_i n_j}{(r_+N_i)(r_j+N_j)}}, & \text{if}\  v_i\sim v_j \ \text{in} \ \mathscr{G}, \\
            0, & \text{otherwise}.
\end{cases}
\end{equation*} 
\end{itemize}
\end{theorem}

\section{Induced Subgraph of $\mathcal{E}_{\mathbb Z_n}$ as a generalized join graph}

In \cite{jamsheena2021structural},  the authors have established the structure of $\mathcal{E}_{\mathbb Z_n}$ as the join of a complete graph and a multipartite graph. 
 In this section, we further investigate this multipartite graph $H$ and prove
that $H$ is the generalized join of certain null graphs.\\

Let the elements of $\mathbb Z_n$ be represented  by  $\{ 0,1,2, \cdots,n\}$. Assume that $n\ne 1$ and $\mathbb Z_n$ is not an integral domain to exclude the triviality of being an empty graph. That is, $n\ge 4$ and $n$ is not a prime. 
If $x\in \mathbb Z_n$, then $\langle x \rangle$ denotes the ideal generated by $x$.
Let  $n= p_{1}^{m_1}p_{2}^{m_2} \cdots p_{k}^{m_k}$, where $p_{1}<p_{2}< \cdots < p_{k}$ be primes, and $m_j$ is a non-negative integer for $1 \le j \le k$. Then, the total number of nonzero proper ideal of $\mathbb Z_n$ is given by $T= \prod_{j=1}^{k}(m_{j}+1)-2$. 
Now, we discuss the structure of $H$ more precisely by defining an equivalence relation on the set $\mathscr{U}$ of nonessential ideals of $\mathbb Z_n$. By Theorem \ref{chr. essentl.idl Zn}, any nonessential ideal of $\mathbb Z_n$ is of the form $I= \langle p_{1}^{r_1}p_{2}^{r_2} \cdots p_{k}^{r_k}\rangle$, where $r_i= m_i$ for at least one $i$ and $r_i\ne m_i$ for all $i$; $1\le i\le k$.
\begin{definition}
     Let $\Xi= \{ 1,2, \cdots, k\}$ be an index set.
    For an ideal $I$ of $\mathscr{U}$, define a subset  $\Xi_I$ of $\Xi$  by,
    $\Xi_I= \{i: r_i= m_i \ \text{in} \ I\}$.\\
    For example, if $I=  \langle p_{1}^{r_1}p_{2}^{r_2} \cdots p_{i-1}^{r_{i-1}} p_{i}^{m_i}p_{i+1}^{r_{i+1}}\cdots p_{t-1}^ {r_{t-1}} p_{t}^{m_t}p_{t+1}^{r_{t+1}} \cdots p_{k}^{r_k} \rangle$,  where $0\le r_s < m_s$ for $1 \le s \le k$ and $s\ne i, t$, then $\Xi_ I= \{i,t \}$.
\end{definition}
    
\begin{definition} \label{eqvlnc reln}
    Let $I$ and $J$ be any two ideals of $\mathscr{U}$. We define a relation $\preccurlyeq$ on $\mathscr{U}$ by $I\preccurlyeq J$ if and only if $\Xi_I= \Xi_J$.
\end{definition}
 We can easily verify that $\preccurlyeq$ is an equivalence relation on $\mathscr{U}$ and it partitions the set  $\mathscr{U}$ of nonessential ideals into nonempty distinct equivalence classes. Since $\Xi$ possesses $2^k -2$ nontrivial proper subsets, there are $2^k -2$ equivalent classes. Let $[I]$ denote the equivalent class of the ideal $I$. \\
 For example, if $I=\langle p_{1}^{m_1}\rangle$ then, $[I]= \{K: K=\langle p_{1}^{m_1}p_{2}^{r_2} \cdots p_{k}^{r_k}\rangle, \ 0\le r_s < m_s \ \text{for}\ 2 \le s \le k\}$. Here, $\Xi_I= \{1\} = \Xi_K$. 
 \begin{lemma}\label{|[I]|}
   $|[I]|=\displaystyle \prod_{i\notin \Xi_I}m_i$, for any equivalent class $[I]$ of the equivalence relation $\preccurlyeq$ on $\mathscr{U}$ .
 \end{lemma}
\begin{lemma}\label{Chrczn Adcncy}
Let $K$ and $L$ be two vertices of any two of the $2^k -2$ equivalent classes, say $[I]$ and $[M]$ respectively. Then $K$ and $L$ are adjacent in $\mathcal{E}_{\mathbb Z_n}$  if and only if $\Xi_I\cap \Xi_M= \phi$.
\end{lemma}
\begin{proof}
   The vertices $K$ and $L$ are of the form $\langle  p_{1}^{r_1}p_{2}^{r_2} \cdots p_{k}^{r_k}\rangle$, where $r_i=m_i$ for $i\in \Xi_I$ and $\langle  p_{1}^{s_1}p_{2}^{s_2} \cdots p_{k}^{s_k}\rangle$, where $s_j=m_j$ for $j\in \Xi_M$ respectively. Then, $K$ and $L$ are adjacent in $\mathcal{E}_{\mathbb Z_n}$ if and only if $K+L=\langle p_{1}^{min\{r_1,s_1\}}p_{2}^{min\{r_2,s_2\}}\cdots\\
   p_{k}^{min\{r_k,s_k\}}\rangle$ is an essential ideal. 
   By Theorem \ref{chr. essentl.idl Zn}, this is possible if and only if $r_t$ and $s_t$ can't be equal to $m_t$ for the same $t$, $1\le t \le k$.
\end{proof}
\begin{corollary}\label{adj of Eqvlt Cls}
\begin{itemize}
\item [(i)] The induced subgraph $\mathcal{E}_{\mathbb Z_n}([I])$ of $\mathcal{E}_{\mathbb Z_n}$ on the vertex set $[I]$ is the complement graph $\overline{K}_{\displaystyle \prod_{i\notin \Xi_I}m_i}$.
 \item [(ii)] For any two equivalent classes $[I]$ and $[M]$, a vertex of $[I]$ is adjacent to either all or none of the vertices of $[M]$ in $\mathcal{E}_{\mathbb Z_n}$. 
\end{itemize}
\end{corollary}
\begin{proof}
   (i) If  $M$ and $L$ are any two vertices of the same equivalent class $[I]$, then they are nonadjacent by Lemma \ref{Chrczn Adcncy}. Also, Lemma \ref{|[I]|} assures that $[I]$ contains $\displaystyle \prod_{i\notin \Xi_I}m_i$ vertices. Hence, $\mathcal{E}_{\mathbb Z_n}([I])\cong \overline{K}_{\displaystyle \prod_{i\notin \Xi_I}m_i}$.\\
   (ii) The proof follows directly from Lemma \ref{Chrczn Adcncy}. 
\end{proof}
Now, Corollary \ref{adj of Eqvlt Cls} assures that the partition of vertices of $\mathscr{U}$ corresponding to the nonessential ideals of $\mathbb Z_n$ into $2^k -2$ equivalent classes are equitable. That is, every vertex in $[I]$ has the same number of neighbors in $[K]$, for any two equivalent classes  $[I]$ and $[K]$. Hence we attempt to prove that  $\mathcal{E}_{\mathbb Z_n}(\mathscr{U})\cong \mathscr{G}[\Gamma_1, \Gamma_2, \cdots, \Gamma_{2^k-2}]$, where each   $\Gamma_i$  represents the subgraph induced by the  equivalent class $[I]$ of the partition resulting from the Definition \ref{eqvlnc reln} on the set $\mathscr{U}$.   \\

\textbf{Construction of the graph $\mathscr{G}$} \\
Let $\mathscr{G}$ be a simple graph where the vertex set comprises representatives from the $2^k - 2$ equivalent classes resulting from the partition of the set $\mathscr{U}$. More precisely, 
\begin{equation}\label{Eq1}
\begin{split}
    V(\mathscr{G})= &\{\langle p_1^{m_1}\rangle,\cdots, \langle p_k^{m_k}\rangle, \langle p_1^{m_1} p_2^{m_2}\rangle,\cdots, \langle p_{k-1}^{m_{k-1}} p_k^{m_k}\rangle, \cdots, \\
    &\langle p_1^{m_1}\cdots p_{k-1}^{m_{k-1}}\rangle,\cdots, \langle p_2^{m_2}\cdots p_k^{m_k}\rangle \},  
\end{split}
\end{equation}
 two vertices of $V(\mathscr{G})$
 are connected by an edge if and only if their corresponding index subsets
have no elements in common.
\begin{lemma}
    The graph $\mathscr{G}$ is connected.
\end{lemma}
\begin{proof}
    Let $I$ and $K$ be any two vertices of  $\mathscr{G}$. If $I\sim K$, then we are done. Suppose $I\nsim K$ in $\mathscr{G}$. Then,  $\Xi_I\cap \Xi_K\ne \phi$. Without loss of generality assume that
    \begin{equation*}
         \Xi_I=\{x_1, x_2,\cdots,x_l,y_1,y_2,\cdots,y_m \}\  \text{and} \
 \Xi_K= \{x_1, x_2,\cdots,x_l,z_1,z_2,\cdots,z_t \}
    \end{equation*} such that $y_i\ne z_j$ for $1\le i\le m$ and $1\le j\le t$. Here, $\Xi_I\cap \Xi_K= \{ x_1, x_2,\cdots,x_l\}$. To show that $I$ and $K$ are connected in $\mathscr{G}$, we consider two cases.\\
    Case $1$: $\Xi \ne \Xi_I \cup \Xi_K$ \\
    Then, $\Xi$ contains at least one element not in $\Xi_I \cup \Xi_K$ which accelerate the existence of at least one vertex $P$ with $\Xi_P\subseteq \Xi\backslash (\Xi_I \cup \Xi_K)$. Hence $I\sim P\sim K$.\\
    Case $2$: $\Xi = \Xi_I \cup \Xi_K$ \\ 
    Since $\langle 0 \rangle$ is not a vertex of $\mathcal{E}_{\mathbb Z_n}([I])$,we have $l+m<k$ and $l+t<k$. So there exists at least one integer, say $e(1\le e\le k)$ such that $e\in \Xi \ \text{but} \ e\notin \Xi_I$. Then, $e$ must belongs to $\Xi_K$. With the same assertion, we claim the existence of an integer $h\ne e$ such that $h\in \Xi \ \text{but}\ h \notin \Xi_K$. Now, take the two vertices $P$ and $Q$ of $\mathscr{G}$ with $\Xi_P= \{e \}$ and $\Xi_Q = \{h\}$ respectively. Hence $I\sim P \sim Q \sim K$.  
\end{proof}
\begin{theorem}\label{G-join1}
    The subgraph of $\mathcal{E}_{\mathbb Z_n}$ induced by the set $\mathscr{U}$ 
    is the generalized join of certain null graphs given by,
\begin{equation*}
\begin{split}
    \mathcal{E}_{\mathbb Z_n}(\mathscr{U}) = \mathscr{G}&[\mathcal{E}_{\mathbb Z_n}([\langle p_1^{m_1}\rangle]),\cdots, \mathcal{E}_{\mathbb Z_n}([\langle p_k^{m_k}\rangle]), \mathcal{E}_{\mathbb Z_n}([\langle p_1^{m_1}p_2^{m_2}\rangle]),\cdots,  \\
    & \mathcal{E}_{\mathbb Z_n}([\langle p_{k-1}^{m_k-1}p_k^{m_k}\rangle]),\cdots, \mathcal{E}_{\mathbb Z_n}([\langle p_2^{m_2}\cdots p_{k-1}^{m_{k-1}}p_k^{m_k} \rangle])],
\end{split}
\end{equation*}
where $\mathcal{E}_{\mathbb Z_n}([I])= \overline{K} _ {\displaystyle \prod_{i\notin \Xi_I}m_i}$ for the representative ideal $I (\text{vertex of}\ \mathscr{G})$ of the equivalent class $[I]$.
\end{theorem}
\begin{proof}
    The proof follows from Lemma \ref{Chrczn Adcncy} by replacing  the vertex $I$ of $\mathscr{G}$ with the induced subgraph $\mathcal{E}_{\mathbb Z_n}([I])$ for each of the $2^k -2$ vertices of $\mathscr{G}$.
\end{proof}
\begin{example}
   Let $n=p_1^{m_1}p_2^{m_2}p_3^{m_3}$, where $p_1<p_2<p_3$ be distinct primes and $m_i$ are positive integers such that $m_i>1$ for at least $1\le i\le 3$. Then the induced subgraph  $\mathcal{E}_{\mathbb Z_n}(\mathscr{U})$ is shown in Figure \ref{fig:1}. Here,\begin{equation*}
   V(\mathscr{G})= \{\langle p_1^{m_1}\rangle,\langle  p_2^{m_2}\rangle , \langle p_3^{m_3}\rangle, \langle p_1^{m_1} p_2^{m_2}\rangle, \langle p_1^{m_1} p_3^{m_3}\rangle, \langle p_2^{m_2} p_3^{m_3}\rangle\} 
 \end{equation*} and $\mathscr{G}$ is the \textit{corona} $K_3\circ K_1$(that is, the graph obtained from $K_3$ by attaching a leaf to every vertex of $K_3$). Also,   
 $\mathcal{E}_{\mathbb Z_n}([\langle p_1^{m_1}\rangle])= \overline{K} _ {m_2m_3}$,   $\mathcal{E}_{\mathbb Z_n}([\langle p_2^{m_2}\rangle])= \overline{K} _ {m_1 m_3}$,   $\mathcal{E}_{\mathbb Z_n}([\langle p_3^{m_3}\rangle])= \overline{K} _ {m_1 m_2}$,  $\mathcal{E}_{\mathbb Z_n}([\langle  p_1^{m_1} p_2^{m_2} \rangle])= \overline{K} _ {m_3}$, $\mathcal{E}_{\mathbb Z_n}([\langle  p_1^{m_1} p_3^{m_3} \rangle])= \overline{K} _ {m_2}$, and $\mathcal{E}_{\mathbb Z_n}([\langle  p_2^{m_2} p_3^{m_3} \rangle])= \overline{K} _ {m_1}$.
\end{example}
\begin{figure}
\begin{minipage}[c]{5 cm}
	\begin{tikzpicture}
\draw[line width=0.25mm,double,color= black!65] (6.72,2.43) -- (7.78,3.71);
		\filldraw [color= blue, opacity=0.20 ] (8,4.15) circle (0.5cm);
		\draw[line width=0.25mm,double,color= black!65] (8.22,3.71) -- (9.28,2.46);
		\draw[line width=0.25mm,double,color= black!65] (8,4.65) -- (8,5.3);
		\filldraw [color= blue, opacity=0.20 ] (8.02,5.8) circle (.5cm);
  \filldraw [color= blue, opacity=0.20](9.5,2.0) circle (.5cm);
		\draw[line width=0.25mm,double,color= black!65] (9.98,2.0) -- (10.64,2);
	\filldraw [color=blue, opacity=0.20 ](11.14,2.0) circle (.5cm);
		\draw[line width=0.25mm,double,color= black!65] (7,2.0) -- (9,2);
		\filldraw [color= blue, opacity=0.20 ] (6.5,2.0) circle (.5cm);
		\draw[line width=0.25mm,double,color= black!65] (6,2) -- (5.36,2);
		\filldraw [color= blue, opacity=0.20 ] (4.85,2.0) circle (.5cm);
		\node at (8,5.8) {\tiny$\overline{K}_{m_1}$};
		\node at (8,4.15) {\tiny$\overline{K}_{m_2m_3}$};
		\node at (9.5,2) {\tiny$\overline{K}_{m_1m_2}$};
		\node at (11.2,2) {\tiny$\overline{K}_{m_3}$};
		\node at (6.5,2) {\tiny$\overline{K}_{m_1m_3}$};
		\node at (4.9,2) {\tiny$\overline{K}_{m_2}$};

\filldraw [color=blue, opacity=0.20 ] (13.5, 2) circle (0.5cm);
\filldraw [color=blue, opacity=0.20 ] (16.5, 2) circle (0.5cm);
\filldraw [color=blue, opacity=0.50 ](14.9, 4) circle (0.5cm);
\draw[line width=0.25mm,double,color= black!65](14,2)--(16,2);
\draw[line width=0.25mm,double,color= red!65](13.6,2.5)--(14.5,3.7);
\draw[line width=0.25mm,double,color= red!65](16.4,2.5)--(15.35,3.89);
\node at (13.5,2) {\tiny$\overline{K}_{m_2}$};
 \node at (16.5,2) {\tiny$\overline{K}_{m_1}$}; 
  \node at (15,4) {\tiny${K}_{m}$}; 
  \end{tikzpicture}
\end{minipage}

\caption{$\mathcal{E}_{\mathbb Z_n}(\mathscr{U})$ for $n= p_1^{m_1}p_2^{m_2}p_3^{m_3}$ \hskip 50pt Figure $2: \mathcal{E}_{\mathbb Z_n}$ for $n= p_1^{m_1}p_2^{m_2}$}\label{fig:1}



\end{figure}


 We further observe that the graph $\mathscr{G}$  of the $\mathscr{G}$-generalized join stated in Theorem \ref{G-join1} has the structure of \textit{Annihilating Ideal Graph} ($\mathbb{AIG}(R)$) associated with a certain commutative ring $R$ found in the literature.
\begin{theorem}\label{G-join2}
   Let  $n= p_{1}^{m_1}p_{2}^{m_2} \cdots p_{k}^{m_k}$ where $p_{1}<p_{2}< \cdots < p_{k}$ are primes, and $m_i$ is a non-negative integer for $1 \le i \le k$ such that $m_i>1$ for at least one $i$. Then, the graph  $\mathscr{G}$ of the $\mathscr{G}$-generalized join of induced subgraph $\mathcal{E}_{\mathbb Z_n}(\mathscr{U})$ is the annihilating ideal graph of the ring $\mathbb Z_n$, where $n$ is the product of $k$ distinct primes. That is, $\mathscr{G}\cong \mathbb{AIG}(\mathbb Z_{n})$ for $n=\displaystyle \prod_{i=1}^{k}p_i$.
\end{theorem}
\begin{proof}
    The annihilating ideal graph $\mathbb{AIG}(\mathbb Z_{n})$, where $n =\prod_{i=1}^{k}p_i$, is a simple connected graph with 
    vertex set \begin{equation*}
    V(\mathbb{AIG}(\mathbb Z_{n}))=\{ I=\langle x \rangle : x \ \text{is a proper divisor of}\ n \}
    \end{equation*} and two vertices $\langle x\rangle$ and $\langle y \rangle$ of  are adjacent if and only if $n \vert xy$.
    Obviously, $\vert  V(\mathbb{AIG}(\mathbb Z_{n})) \vert= 2^k-2$. To prove $\mathscr{G}\cong \mathbb{AIG}(\mathbb Z_{n})$, consider the mapping $\psi: V(\mathscr{G})\rightarrow  V(\mathbb{AIG}(\mathbb Z_{n}))$ defined by $\psi(I)= \langle \displaystyle \prod_{i\notin \Xi_I} p_i \rangle$. We can easily verify that $\psi$ is well-defined. \\
    \textbf{claim}:  $\psi$ is an isomorphism\\
    At first, we prove that $\psi$ is one-one and onto. Consider any two distinct vertices $I$ and $K$ of $ V(\mathscr{G})$ for this.
    Then, $\Xi_I\ne \Xi_K$ leading to two cases: either  $ \Xi_I\cap \Xi_K = \phi$ or $ \Xi_I\cap \Xi_K \ne \phi$.\\
    Case $1$: $ \Xi_I\cap \Xi_K = \phi$\\
    If $\Xi_I\cup \Xi_K= \Xi$, then $\psi(I)$ and $\psi(K)$ are distinct.
    If not, $\psi(I)$ and $\psi(K)$ are given by, 
    \begin{equation*}
        \psi(I)= \langle  \prod_{i\in \Xi_K}p_i \prod_{i\in  \Xi\backslash (\Xi_I \cup \Xi_K)}p_i \rangle \ \text{and}\ 
         \psi(K)= \langle  \prod_{i\in \Xi_I}p_i \prod_{i\in  \Xi\backslash (\Xi_I \cup \Xi_K)}p_i \rangle
    \end{equation*}respectively. There is at least one distinct prime in the products $\displaystyle\prod_{i\in \Xi_I}p_i$ and $\displaystyle \prod_{i\in \Xi_K}p_i$ which guarantees that $\psi(I)$ and $\psi(K)$ are distinct. \\
    Case $2$:  $ \Xi_I\cap \Xi_K \ne  \phi$\\
    In this case, if $\Xi_I\cup \Xi_K= \Xi$, we can see that the ideals  $\psi(I)$ and $\psi(K)$ are distinct. If  $\Xi_I\cup \Xi_K\ne \Xi$,  $\psi(I)$ and $\psi(K)$ given by  
    \begin{equation*}
        \psi(I)=\langle \prod_{i\in \Xi_{K}\backslash (\Xi_I \cap \Xi_K)}p_i\prod_{i\in  \Xi\backslash (\Xi_I \cup \Xi_K)}p_i \rangle
    \end{equation*} and 
    \begin{equation*}
         \psi(K)= \langle  \prod_{i\in \Xi_{I}\backslash (\Xi_I \cap \Xi_K)}p_i  \prod_{i\in  \Xi\backslash (\Xi_I \cup \Xi_K)}p_i \rangle
    \end{equation*} are  distinct. \\
Also, for any vertex $K'=\langle \displaystyle \prod_{i\in \Xi_{K'}}p_i\rangle$  of $V(\mathbb{AIG}(\mathbb Z_{n}))$, $\Xi_{K'}\subset \Xi$. Hence, there exists an ideal $K\in V(\mathscr{G})$ with index subset $\Xi_{K}= {\Xi_{K'}}^c$ such that $\psi(K)= \langle \displaystyle \prod_{i\notin \Xi_{K}} p_i\rangle= K'$. Thus $\psi$ is one-one and onto. \\
\textit{$\psi$ preserves adjacency}: Let $I$  and $K$ be two adjacent vertices in $\mathscr{G}$  with index sets $\Xi_I$ and $\Xi_K$ respectively. Then, $ \Xi_I\cap \Xi_K = \phi$.
If $\Xi_I\cup \Xi_K= \Xi$, the vertices $\psi(I)$ and $\psi(K)$ are adjacent in $\mathbb{AIG} (\mathbb Z_{n})$, since $\psi(I)\psi(K)= \langle 0 \rangle$. Also, if  $\Xi_I\cup \Xi_K\ne \Xi$, 
\begin{equation*}
    \psi(I)\psi(K)= \langle \prod_{i\in \Xi}p_i \prod_{i\in  \Xi\backslash (\Xi_I \cup \Xi_K)}p_i \rangle=\langle 0\rangle, 
\end{equation*} proves that $\psi(I)$ and $\psi(K)$ are adjacent in $\mathbb{AIG}(\mathbb Z_{n})$.        

\end{proof}
As a consequence of Theorems \ref{Strctr thm Zn}, \ref{G-join1} and \ref{G-join2}, the structure of  $\mathcal{E}_{\mathbb Z_n}$ is completely determined in the following theorem.
\begin{theorem} \label{Complt Strctr of EZn}
   Let  $n= p_{1}^{m_1}p_{2}^{m_2} \cdots p_{k}^{m_k}$, where $p_{1}<p_{2}< \cdots < p_{k}$ are primes, and $m_i$ is a non-negative integer for $1 \le i \le k$ such that $m_i>1$ for at least one $i$. Then, the essential ideal graph  $\mathcal{E}_{\mathbb Z_n}\cong K_m \vee H$, where $K_m$ is the complete graph  on $m=\prod_{i=1}^{k}m_{i}-1$ vertices and  \begin{equation*}
       \begin{split}
      H= \mathscr{G}&[\mathcal{E}_{\mathbb Z_n}([\langle p_1^{m_1}\rangle]),\cdots, \mathcal{E}_{\mathbb Z_n}([\langle p_k^{m_k}\rangle]), \mathcal{E}_{\mathbb Z_n}([\langle p_1^{m_1}p_2^{m_2}\rangle]),\cdots,  \\
    & \mathcal{E}_{\mathbb Z_n}([\langle p_{k-1}^{m_k-1}p_k^{m_k}\rangle]),\cdots,\mathcal{E}_{\mathbb Z_n}([\langle p_2^{m_2}\cdots p_{k-1}^{m_{k-1}}p_k^{m_k} \rangle])];
\end{split}
\end{equation*} $ \mathscr{G}\cong \mathbb{AIG}(\mathbb Z_{\prod_{i=1}^{k}p_i})$.
\end{theorem}
    \section{Spectra of the induced subgraph $\mathcal{E}_{\mathbb Z_n}(\mathscr{U})$}   

In this section, we determine the various spectrum such as adjacency spectrum, Laplacian spectrum, and signless Laplacian spectrum of the induced subgraph  $\mathcal{E}_{\mathbb Z_n}(\mathscr{U})$. Also, a characterization for $\mathcal{E}_{\mathbb Z_n}$ to be  Laplacian integral is obtained.
\begin{theorem}\label{Aspec EZn}
 Let  $n= p_{1}^{m_1}p_{2}^{m_2} \cdots p_{k}^{m_k}$, where $p_{1}<p_{2}< \cdots < p_{k}$ are primes, and $m_i$ is a non-negative integer  such that $m_i>1$ for at least one $i$ . Then the eigenvalues of $\mathcal{E}_{\mathbb Z_n}(\mathscr{U})$ are $0$ with multiplicity $T-m-(2^k-2)$,
 and the remaining $2^k-2$ eigenvalues are contained in the spectrum of the matrix,\\
 $C_A(\mathscr{G})=(c_{ij})= \begin{cases}
       0, & \text{if}\ i=j,\\
       \sqrt{n_In_J}, &  \text{if}\ I\sim J \ \text{in}\ \mathscr{G},\\
       0, &\text{otherwise}.\\
   \end{cases}$ 
\end{theorem}
\begin{proof}
    By Theorem \ref{G-join1}, \begin{equation*}
\begin{split}
    \mathcal{E}_{\mathbb Z_n}(\mathscr{U}) = \mathscr{G}&[\mathcal{E}_{\mathbb Z_n}([\langle p_1^{m_1}\rangle]),\cdots, \mathcal{E}_{\mathbb Z_n}([\langle p_k^{m_k}\rangle]), \mathcal{E}_{\mathbb Z_n}([\langle p_1^{m_1}p_2^{m_2}\rangle]),\cdots,  \\
    & \mathcal{E}_{\mathbb Z_n}([\langle p_{k-1}^{m_k-1}p_k^{m_k}\rangle]),\cdots, \mathcal{E}_{\mathbb Z_n}([\langle p_2^{m_2}\cdots p_{k-1}^{m_{k-1}}p_k^{m_k} \rangle])],
\end{split}
\end{equation*} where $\mathcal{E}_{\mathbb Z_n}([I])= \overline{K} _ {n_I}$ 
 Also, note that $\sigma_A( \overline{K} _ {n_I})=\begin{pmatrix}
    0 \\
    n_I \\
\end{pmatrix}$. Then the result follows immediately from Theorem \ref{spec G-join}.  
\end{proof}

\begin{theorem}\label{LSN spec EZn}
    Let  $n= p_{1}^{m_1}p_{2}^{m_2} \cdots p_{k}^{m_k}$, where $p_{1}<p_{2}< \cdots < p_{k}$ are primes, $m_i$ is a non-negative integer  such that $m_i>1$ for at least one $i$, $n_I=|\mathcal{E}_{\mathbb Z_n}([I])|$, and let $N_I=\displaystyle \sum_{J}n_J$ where the summation runs over all the ideals $J\in V(\mathscr{G})$ for which  $\Xi_I\cap \Xi_J= \phi$. Then,
    \begin{enumerate}
        \item \label{LSpec EZn} the Laplacian eigenvalues of $\mathcal{E}_{\mathbb Z_n}(\mathscr{U})$ are  $N_I$ with multiplicity \\
      $\bigg(\displaystyle \prod_{i\notin \Xi_I}m_i-1 \bigg),\  I\in V(\mathscr{G})$
  for which $\mathcal{E}_{\mathbb Z_n}([I])$
       is not  singleton  and the remaining are the eigenvalues of the matrix 
    \begin{equation*}
            C_L(\mathscr{G})=(c_{ij})= \begin{cases}
             N_I, & \text{if}\ i=j, \\
              -\sqrt{n_In_J},  & \text{if}\ I\sim J \text{in}\ \mathscr{G}, \\
              0, & \text{otherwise}.
    \end{cases}
    \end{equation*}
        \item \label{SLSpec EZn}the signless Laplacian eigenvalues of  $\mathcal{E}_{\mathbb Z_n}(\mathscr{U})$ contains the same  \[ \text{Laplacian eigenvalues}\ N_I\ \text{with multiplicity}\ 
    \bigg(\prod_{i\notin \Xi_I}m_i-1 \bigg),\  I\in V(\mathscr{G}),\] and the remaining are eigenvalues of the matrix 
       \begin{equation*}
            C_Q(\mathscr{G})=(c_{ij})= \begin{cases}
             N_I, & \text{if}\ i=j, \\
              \sqrt{n_In_J},  & \text{if}\ I\sim J \text{in}\ \mathscr{G}, \\
              0, & \text{otherwise}.
    \end{cases}
    \end{equation*}
    \item \label{NLSpec EZn}the normalized Laplacian eigenvalues of  $\mathcal{E}_{\mathbb Z_n}(\mathscr{U})$ are,
    \[1 \ \text{with multiplicity} \sum_{I}(\prod_{i\notin \Xi_I}m_i -1)\ \text{in which the summation runs over} \] all the ideals $I\in V(\mathscr{G})$, and the remaining are the eigenvalues of the matrix
    \begin{equation*}
        C_{\mathscr{L}}(\mathscr{G})= (c_{ij})= \begin{cases}
          1, & \text{if}\ i=j, \\
              -\sqrt{\frac{n_In_J}{N_IN_J}},  & \text{if}\ I\sim J \text{in}\ \mathscr{G}, \\
              0, & \text{otherwise}.   
        \end{cases}
    \end{equation*}
    \end{enumerate}
\end{theorem}
\begin{proof}
The proof follows from Theorems  \ref{spec G-join} and  \ref{G-join1}.
\end{proof}
\begin{example}
 Let $n=p^2qr$. Then $V(\mathscr{G})= \{\langle p^2\rangle, \langle q \rangle, \langle r \rangle, \langle p^2q \rangle, \langle p^2r \rangle, \langle qr \rangle \}$, $\mathcal{E}_{\mathbb Z_n}([\langle p^2\rangle])= \overline{K} _ 1=\mathcal{E}_{\mathbb Z_n}([\langle p^2q\rangle])=\mathcal{E}_{\mathbb Z_n}([\langle p^2r\rangle])$,\\
 $\mathcal{E}_{\mathbb Z_n}([\langle q\rangle])= \overline{K}_2=\mathcal{E}_{\mathbb Z_n}([\langle r\rangle])= \mathcal{E}_{\mathbb Z_n}([\langle qr\rangle])$, and \\
    $n_{\langle p^2\rangle}= 1$, $n_{\langle q \rangle}=2$, $n_{\langle r \rangle}=2$,
    $n_{\langle p^2q\rangle}= 1$, $n_{\langle p^2r \rangle}= 1$ and $n_{\langle qr \rangle}=2$. \\
    $N_{\langle p^2 \rangle}= n_{\langle q \rangle} + n_ {\langle r \rangle} + n_ {\langle qr \rangle}= 6$,  $N_{\langle q \rangle}= n_{\langle p^2 \rangle} + n_ {\langle r \rangle} + n_ {\langle p^2r \rangle}= 4= N_{\langle r \rangle} $,\\
     $N_{\langle p^2q \rangle}= n_ {\langle r \rangle}= 2= N_{\langle p^2r \rangle}$, and $N_{\langle qr \rangle}= n_ {\langle p^2 \rangle}= 1$.
    \begin{itemize}
 \item By Theorem \ref{Aspec EZn}, $0$ is an eigenvalue of multiplicity $3$ and the remaining eigenvalues are that of the $6\times 6$ matrix 
    \[ C_A(\mathscr{G})=  \left(
\begin{array}{c|c}
    \sqrt{2}(J-I)_{3\times3} & B_{3\times3} \\ \hline
B_{3\times3} & 0_{3\times3} \\ 
\end{array}
\right),\] where $B= \begin{pmatrix}
    0 & 0 & \sqrt{2}\\
     0 & \sqrt{2}& 0 \\
    \sqrt{2}& 0 & 0 \\
\end{pmatrix}$. By Lemma \ref{detMATX}, we obtain the spectrum of $C_A(\mathscr{G})$ as,
\begin{equation*}
    \sigma(C_A(\mathscr{G}))=\begin{pmatrix}
        2+\sqrt{2} & \frac{-1+\sqrt{1+2\sqrt{2}}}{\sqrt{2}}& -(2-\sqrt{2}) &  \frac{-1-\sqrt{1+2\sqrt{2}}}{\sqrt{2}}\\
        1 &2 & 1& 2 \\
    \end{pmatrix},
\end{equation*} and hence
\begin{equation*}
    \sigma_A(\mathcal{E}_{\mathbb Z_n}(\mathscr{U}))= \begin{pmatrix}
        2+\sqrt{2} & \frac{-1+\sqrt{1+2\sqrt{2}}}{\sqrt{2}} & 0 & -(2-\sqrt{2}) &  \frac{-1-\sqrt{1+2\sqrt{2}}}{\sqrt{2}}\\
        1 & 2 & 3 & 1 & 2 \\
    \end{pmatrix}.
\end{equation*}
         \item By Theorem \ref{LSN spec EZn}(\ref{LSpec EZn}), the Laplacian eigenvalues of $\mathcal{E}_{\mathbb Z_n}(\mathscr{U})$ are $4$ with 
 multiplicity $2$, $1$ as a simple eigenvalue
 and the remaining six eigenvalues are of the matrix 
\[ C_L(\mathscr{G})=  \left(
\begin{array}{ccc|cccc}
 6 & -\sqrt{2} &  -\sqrt{2} & 0 & 0 & -\sqrt{2} \\
 -\sqrt{2} & 4 & -2 & 0 & -\sqrt{2} & 0 \\
 -\sqrt{2} & -2 & 4 &  -\sqrt{2} & 0 & 0 \\
 \hline
  0 & 0 & -\sqrt{2} & 2 & 0 & 0 \\
 0 & -\sqrt{2}  & 0 & 0 & 2 & 0\\
-\sqrt{2}  & 0 & 0   & 0 & 0 & 1 \\
\end{array}
\right)=  \left(
\begin{array}{c|c}
C & N \\ \hline
N & D
\end{array}
\right).
\]
For $\lambda\ne 1,2$, by Lemma \ref{detMATX}, we have 
\begin{equation*}
\begin{split}
det( C_L(\mathscr{G})- \lambda I)& =det(D-\lambda I). det(C-ND^{-1}N^T)\\
& = \lambda(\lambda-3)(\lambda^2-8\lambda+6) (\lambda^2-8\lambda+10).    
\end{split}
\end{equation*}
Hence,
\begin{equation*}
    \sigma_L(\mathcal{E}_{\mathbb Z_n}(\mathscr{U}))= \begin{pmatrix}
        4+\sqrt{10} & 4+\sqrt{6} & 4 & 3 &  4-\sqrt{6} & 1 & 4-\sqrt{10} & 0 \\  
         1 & 1 & 2 & 1 & 1 & 1 & 1 & 1 \\
    \end{pmatrix}.
\end{equation*}
\item Using Theorem \ref{LSN spec EZn}(\ref{SLSpec EZn}), the signless Laplacian spectrum of  $\mathcal{E}_{\mathbb Z_n}(\mathscr{U})$ contains the same eigenvalues $4$ with multiplicity $2$, $1$ as a simple eigenvalue of the Laplacian matrix of $\mathcal{E}_{\mathbb Z_n}(\mathscr{U})$, and the remaining six eigenvalues are that of the matrix 
\[ C_Q(\mathscr{G})=  \left(
\begin{array}{ccc|cccc}
 6 & \sqrt{2} &  \sqrt{2} & 0 & 0 & \sqrt{2} \\
 \sqrt{2} & 4 & 2 & 0 & \sqrt{2} & 0 \\
 \sqrt{2} & 2 & 4 &  \sqrt{2} & 0 & 0 \\
 \hline
  0 & 0 & \sqrt{2} & 2 & 0 & 0 \\
 0 & \sqrt{2}  & 0 & 0 & 2 & 0\\
\sqrt{2}  & 0 & 0   & 0 & 0 & 1 \\
\end{array}
\right)\]

Then, 
\begin{equation*}
    det( C_Q(\mathscr{G})- \lambda I)= (\lambda^2-4\lambda+2)(\lambda^4-15\lambda^3+66\lambda^2-90\lambda+32).
\end{equation*}
Hence, 
\begin{equation*}
    \sigma_Q(\mathcal{E}_{\mathbb Z_n}(\mathscr{U}))= \begin{pmatrix}
        8.296 & 4.649 & 4 & 2+\sqrt{2} & 1.503 & 1 & 2-\sqrt{2} & 0.552  \\  
         1 & 1 & 2 & 1 & 1 & 1 & 1 & 1 \\
    \end{pmatrix}.
\end{equation*}
\item  By Theorem \ref{LSN spec EZn}(\ref{NLSpec EZn}), the normalized Laplacian eigenvalues of $\mathcal{E}_{\mathbb Z_n}(\mathscr{U})$ are $1$ having multiplicity $3$ and the remaining six eigenvalues are contained in the spectrum of the matrix  
\[ C_{\mathscr{L}}(\mathscr{G})=  \left(
\begin{array}{ccc|cccc}
1 & -\frac{1}{\sqrt{12}} &  -\frac{1}{\sqrt{12}} & 0 & 0 & -\frac{1}{\sqrt{3}} \\
-\frac{1}{\sqrt{12}} & 1 & -\frac{1}{2} & 0 &  -\frac{1}{2}  & 0 \\
 \sqrt{2} & 2 & 1 &   -\frac{1}{2} & 0 & 0 \\
 \hline
  0 & 0 & -\frac{1}{2} & 1 & 0 & 0 \\
 0 &  -\frac{1}{2} & 0 & 0 & 1 & 0\\
-\frac{1}{\sqrt{3}} & 0 & 0   & 0 & 0 & 1 \\
\end{array}
\right)
\]
Here, \begin{equation*}
    det( C_{\mathscr{L}(\mathscr{G})}- \lambda I)= \lambda(4\lambda^2-10\lambda+5)(12\lambda^3-42\lambda^2+45\lambda-14),
\end{equation*} and hence 
\begin{equation*}
    \sigma_{\mathscr{L}}(\mathcal{E}_{\mathbb Z_n}(\mathscr{U}))= \begin{pmatrix}
        \frac{5+\sqrt{5}}{4} & 1.687 &  1.267 & 1 & \frac{5-\sqrt{5}}{4} & 0.545 & 0 \\  
         1 & 1 & 1 & 3 & 1 & 1 & 1  \\
    \end{pmatrix}.
\end{equation*}
 \end{itemize} 
\end{example}
Next, we discuss a characterization of the graph $\mathcal{E}_{\mathbb Z_n}$  to be  Laplacian integral  in terms of the vertex-weighted Laplacian matrix.\\
Consider  $\mathscr{G}$, the underlying graph of the generalized join of $\mathcal{E}_{\mathbb Z_n}(\mathscr{U})$ defined in Section $3$, as a vertex-weighted graph obtained by assigning the weight $n_I= |\mathcal{E}_{\mathbb Z_n}([I])|$ to each vertex $I$ of $2^k-2$ vertices of $\mathscr{G}$. \\
Let $\textbf{L}(\mathscr{G})= (l_{i,j})$ be a matrix of order $k$, where \begin{equation*}
 l_{i,j}= \begin{cases}
\sum_{I \sim J}n_J, & \text{if} \ i=j, \\
-n_J, &  \text{if} \ i\ne j \ \text{and}\ I \sim J\ \text{in} \ \mathscr{G}, \\
0, & otherwise. \\
\end{cases}    
\end{equation*}
The matrix  $\textbf{L}(\mathscr{G})$ is called a \textit{vertex weighted Laplacian matrix of $\mathscr{G}$}, which is a zero-row sum matrix but not symmetric in general. Even though the matrix $C_L(\mathscr{G})$ of order $2^k-2$ defined in Theorem \ref{LSN spec EZn}(\ref{LSpec EZn}) is symmetric, it need not be a zero-row sum matrix. If $\textit{W}= diag[n_1, n_2,\cdots,n_{2^k-2}]$, then $\textbf{L}(\mathscr{G})= \textit{W}^{-\frac{1}{2}}C_L(\mathscr{G}) \textit{W}^{\frac{1}{2}}$ \cite{chung1996combinatorial}. 
So the matrices $\textbf{L}(\mathscr{G})$ and $C_L(\mathscr{G})$ are similar. \\
We thus have the following result.

\begin{proposition} \label{Eql spec}
    $\sigma(C_L(\mathscr{G})= \sigma(\textbf{L}(\mathscr{G})$.
\end{proposition}
By Theorems \ref{G-join1}, \ref{LSN spec EZn}(\ref{LSpec EZn}) and Proposition \ref{Eql spec}, the Laplacian spectrum of $\mathcal{E}_{\mathbb Z_n}(\mathscr{U})$ can be described as follows.
\begin{theorem}
   Let  $n= p_{1}^{m_1}p_{2}^{m_2} \cdots p_{k}^{m_k}$ where $p_{1}<p_{2}< \cdots < p_{k}$ are primes, $m_i$ is a non-negative integer for $1 \le i \le k$. Then the Laplacian eigenvalues of $\mathcal{E}_{\mathbb Z_n}(\mathscr{U})$ are  \[ N_I\ \text{with multiplicity}\ 
       \bigg(\prod_{i\notin \Xi_I}m_i-1 \bigg), \  I\in V(\mathscr{G})\] for which $\mathcal{E}_{\mathbb Z_n}([I])$ is not singleton and the remaining are eigenvalues of the vertex weighted Laplacian matrix  $\textbf{L}(\mathscr{G})$.
\end{theorem}
\begin{proposition}\label{Lplcn Intgrblty EZn}
    The essential ideal graph $\mathcal{E}_{\mathbb Z_n}$ is Laplacian integral if and only if all the eigenvalues of $\textbf{L}(\mathscr{G})$ are integers.
\end{proposition}
\begin{corollary}
    If $n=p_{1}^{m_1}p_{2}^{m_2}$, $p_1< p_2$ are primes and $m_i>1$ for at least one $i$, then the graph $\mathcal{E}_{\mathbb Z_n}$ is Laplacian integral.
\end{corollary}
\begin{proof}
    By Theorem \ref{Complt Strctr of EZn},  $\mathcal{E}_{\mathbb Z_n}\cong K_m \vee H$, where $m=\bigg(\displaystyle\prod_{i=1}^{2} m_i-1\bigg)$ and $H= \mathscr{G}[\mathcal{E}_{\mathbb Z_n}([\langle p_1^{m_1}\rangle]),\ \mathcal{E}_{\mathbb Z_n}([\langle p_2^{m_2}\rangle]); \mathscr{G}= \mathbb{AIG}(\mathbb Z_{p_1p_2})= K_2$ (see Figure $2$). 
    Here, $n_{\langle p_1^{m_1} \rangle}= m_2,\ n_{\langle p_2^{m_2} \rangle}= m_1,\  N_{\langle p_1^{m_1} \rangle}=n_{\langle p_2^{m_2} \rangle}= m_1,\  N_{\langle p_2^{m_2} \rangle}=m_2$ and 
    \[ \textbf{L}(\mathscr{G})= \begin{pmatrix}
       m_1 & -m_2 \\
       -m_1 & m_2 \\
    \end{pmatrix}.\  \text{Obviously, both eigenvalues of}\ \textbf{L}(\mathscr{G}) \ \text{are integers}.\] Then, the result follows from Proposition \ref{Lplcn Intgrblty EZn}.
\end{proof}

\section{\textbf{Spectral Radius and Algebraic Connectivity of  $\mathcal{E}_{\mathbb Z_n}$}}
In this section, we find the bounds for the algebraic connectivity and spectral radius of $\mathcal{E}_{\mathbb Z_n}$. For this, we make use of the following results from the literature.
\begin{theorem}\label{upbd of b(G)} \cite{fiedler1973algebraic}
    Let $\Gamma$ be a graph on $n$ vertices. Then $b(\Gamma)\le n$. Further, equality holds if and only $\overline{\Gamma}$ is disconnected.
\end{theorem}
\begin{theorem}\label{alg cnncvty=v.cnncvty}\cite{kirkland2002graphs}
    For a non-complete connected graph $\Gamma$ with $n$ vertices, $a(\Gamma)=\kappa(\Gamma)$ if and only if $\Gamma=\Gamma_1 \vee \Gamma_2$, where $\Gamma_1$ is a disconnected graph on $n-\kappa(\Gamma)$ vertices and $\Gamma_2$ is a graph on $\kappa(\Gamma)$ vertices with $a(\Gamma_2) \ge 2 \kappa(\Gamma)-n$. 
\end{theorem}
\begin{lemma}\label{cnctvty of EZn comple}
   Let  $n= p_{1}^{m_1}p_{2}^{m_2} \cdots p_{k}^{m_k}$ where $p_{1}<p_{2}< \cdots < p_{k}$ are primes, and $m_i$ is a non-negative integer for $1 \le i \le k$. Then,   $\overline{\mathcal{E}_{\mathbb Z_n}}$ is connected if and only if $n= p_1p_2\cdots p_k$; $k\ge 3$.   
\end{lemma}
\begin{proof}
Firstly, let $n=p_1p_2\cdots p_k$; $k\ge 3$. Then the vertices of the form $I=\langle p_1p_2\cdots p_{i-1}p_{i+1}\cdots p_k\rangle$ are pendant in $\mathcal{E}_{\mathbb Z_n}$ and $N_{\mathcal{E}_{\mathbb Z_n}}(I)= \{\langle p_i \rangle\}$, for $i\in \{1,2,\cdots, k\}$. Then, $deg(I)=T-2$ in  $\overline{\mathcal{E}_{\mathbb Z_n}}$. Consequently, in  $\overline{\mathcal{E}_{\mathbb Z_n}}$, any two of these $T-2$ vertices are connected by a path between them. Now, to prove that $\overline{\mathcal{E}_{\mathbb Z_n}}$ is connected, the only thing remaining is to show that the two vertices $I$ and $\langle p_i \rangle$ are also connected in $\overline{\mathcal{E}_{\mathbb Z_n}}$. Obviously, for a fixed $i$, the vertex $I$ is adjacent to each of the vertices $J=\langle p_1p_2\cdots p_{j-1}p_{j+1}\cdots p_k\rangle$ for $j\ne i \in \{ 1,2,\cdots,\ k\}$ in  $\overline{\mathcal{E}_{\mathbb Z_n}}$. Hence  $\langle p_i \rangle \sim J \sim I$.

For  $n=p^t, t>2$ and $n=p_1p_2$, $\mathcal{E}_{\mathbb Z_n}= K_{t-1}\  \text{and}\ K_2$ respectively. Thus in both cases, $\overline{\mathcal{E}_{\mathbb Z_n}}$ is disconnected. 
Now, the only case remaining is $n= p_{1}^{m_1}p_{2}^{m_2} \cdots p_{k}^{m_k}$, where $k\ge 2$ and $m_i>1$ for at least one $i$. Here, $\mathcal{E}_{\mathbb Z_n}$ contains at least one essential ideal which will be isolated in $\overline{\mathcal{E}_{\mathbb Z_n}}$.
\end{proof}
\begin{theorem} Let $n= p_{1}^{m_1}p_{2}^{m_2} \cdots p_{k}^{m_k}$ where $p_{1}<p_{2}< \cdots < p_{k}$ are primes, and $m_i$ is a non-negative integer for $1 \le i \le k$.  
  For the graph $\mathcal{E}_{\mathbb Z_n}$, $b(\mathcal{E}_{\mathbb Z_n})\le T$. Here, equality holds if and only if  one of the following holds:
  \begin{itemize}
      \item[(i)] $k=1 \ \text{and}\ m_1=t>2$,
      \item [(ii)]  $k=2$ with $m_1=m_2=1$,
      \item [(iii)]  $k\ge 2$ and $m_i>1$ for at least one $i$.
  \end{itemize}
\end{theorem}
\begin{proof}
    By Theorem \ref{upbd of b(G)}, $b(\mathcal{E}_{\mathbb Z_n})\le T$ and $b(\mathcal{E}_{\mathbb Z_n})=T$ if and only if $\overline{\mathcal{E}_{\mathbb Z_n}}$ is disconnected.  
    Now, Lemma \ref{cnctvty of EZn comple} assures that  $\overline{\mathcal{E}_{\mathbb Z_n}}$ is disconnected whenever $n=\{ p_{1}^t,\  p_1p_2, \ p_{1}^{m_1}p_{2}^{m_2} \cdots p_{k}^{m_k}\}$ for $t>2$, $k\ge 2$ and $m_i>1$ for at least one $i$ respectively.
\end{proof}

The following proposition determines the vertex connectivity of the essential ideal graph of the ring $\mathbb Z_{n}$.
\begin{proposition}\label{alg cnncty of EZn}
  Let $n= p_{1}^{m_1}p_{2}^{m_2} \cdots p_{k}^{m_k}$ where $p_{1}<p_{2}< \cdots < p_{k}$ are primes and $T=|V(\mathcal{E}_{\mathbb Z_n})|$. Then, \begin{equation*}
        \kappa(\mathcal{E}_{\mathbb Z_n})= \begin{cases}
          T-1, &  \text{if}\ k=1 \ \text{and} \ m_1>2 \\
           1, & \text{if}\  k>1\ \text{and}\ m_i=1 \ \text{for all}\ i\\
        m+\eta, &   \text{if}\ k>1 \ \text{and}\  m_i>1 \ \text{for at least one}\  i   \\
\end{cases},
\end{equation*}
\[
\eta= min_{I\in V(\mathscr{G})}\{|[I]|;I=\langle p_{j}^{m_j}\rangle \ \text{for}\ 1\le j\le k \}\].    
\end{proposition}
 \begin{proof}
     Case $1$: $n=p^t; \ t>2$\\
By Theorem \ref{chr. essentl.idl Zn}, $\mathcal{E}_{\mathbb Z_n}= K_{t-1}$. Thus, $\kappa(\mathcal{E}_{\mathbb Z_n})=t-2=T-1$.\\
     Case $2$: $k>1\ \text{and}\ m_i=1 \ \text{for all}\ i$\\
     That is, $n$ is the product of $k$ distinct primes. Since the vertices of the form $I= \langle p_1 p_2\cdots p_{i-1}p_{i+1}\cdots p_k \rangle$ are pendant in $\mathcal{E}_{\mathbb Z_n}$,  all the prime ideals of $\mathbb Z_n$ ($\langle p_i \rangle$) are cut vertices in $\mathcal{E}_{\mathbb Z_n}$ for $1\le i \le k$. \\ 
     Case $3$: By Theorems  \ref{Strctr thm Zn} and \ref{G-join1}, $\mathcal{E}_{\mathbb{Z}_{n}}$ is the join of two graphs $K_m$ and $H$, where  
\[
\begin{split}
   H= \mathscr{G}&[\mathcal{E}_{\mathbb Z_n}([\langle p_1^{m_1}\rangle]),\cdots, \mathcal{E}_{\mathbb Z_n}([\langle p_k^{m_k}\rangle]), \mathcal{E}_{\mathbb Z_n}([\langle p_1^{m_1}p_2^{m_2}\rangle]),\cdots,\\
    & \cdots, \mathcal{E}_{\mathbb Z_n}([\langle p_2^{m_2}\cdots p_{k-1}^{m_{k-1}}p_k^{m_k} \rangle])]\ \text{and}\  m=(\prod_{i=1}^{k}m_i - 1).
\end{split} 
\] In the generalized join graph $H$, the vertices of the induced subgraph $\mathcal{E}_{\mathbb Z_n}([I])$ for
$I= \langle p_1^{m_1}\cdots p_{i-1}^{m_{i-1}} p_{i+1}^{m_{i+1}}\cdots p_{k-1}^{m_{k-1}}p_k^{m_k} \rangle $ are adjacent only to the vertices of the induced subgraph  $\mathcal{E}_{\mathbb Z_n}([\langle p_{i}^{m_i}\rangle])$, $1\le i \le k$.
Then, \[\kappa(H)= min \{|[\langle p_{i}^{m_i}\rangle]|,\ 1\le i \le k\}= \eta. \] Now, 
\[\begin{split}
   \kappa(\mathcal{E}_{\mathbb{Z}_{n}}) & =min\{m+\kappa(H), T-m+\kappa(K_m)\} \\
   & = min\{m+\eta, T-1\}= m+\eta. \\
\end{split} \]
\end{proof}
 For non-complete graphs, the algebraic connectivity $a(\Gamma)\le \kappa(\Gamma)$ \cite{fiedler1973algebraic}. Obviously, $a(\mathcal{E}_{\mathbb{Z}_{n}})\le \kappa(\mathcal{E}_{\mathbb{Z}_{n}})$ for $n \ne p^t,\ t>2$ and $n\ne pq$, where $p<q$ be distinct primes. Now, we determine the cases for which $a(\mathcal{E}_{\mathbb{Z}_{n}})< \kappa(\mathcal{E}_{\mathbb{Z}_{n}})$ and $a(\mathcal{E}_{\mathbb{Z}_{n}})= \kappa(\mathcal{E}_{\mathbb{Z}_{n}})$. 
\begin{proposition}
\begin{itemize}
        \item[1]. If $n= p_1p_2 \cdots p_k$; $k>2$, then $a(\mathcal{E}_{\mathbb{Z}_{n}})< \kappa(\mathcal{E}_{\mathbb{Z}_{n}})$. 
        \item [2]. If $n=p_1^{m_1}p_2^{m_2}$, where $m_i>1$ for at least $1\le i\le 2$, then $a(\mathcal{E}_{\mathbb{Z}_{n}})= \kappa(\mathcal{E}_{\mathbb{Z}_{n}})$. 
\end{itemize}
\end{proposition}
\begin{proof}

\begin{enumerate}
    \item If   $n= p_1p_2 \cdots p_k$, then $\overline{\mathcal{E}_{\mathbb Z_n}}$ is connected by Lemma \ref{cnctvty of EZn comple}. So $\mathcal{E}_{\mathbb{Z}_{n}}$ can't be a join of two graphs. Also, $\mathcal{E}_{\mathbb{Z}_{n}}$ is non-complete and connected. Hence by Theorem \ref{alg cnncvty=v.cnncvty}, $a(\mathcal{E}_{\mathbb{Z}_{n}})< \kappa(\mathcal{E}_{\mathbb{Z}_{n}})$. 
\item In this case, $\mathcal{E}_{\mathbb{Z}_{n}}\cong K_m \vee H$, where $H= \overline{K_{m_1}}\vee \overline{K_{m_2}}$. By  Theorem \ref{LcharPlnml join}, 
\begin{equation*}\begin{split}
    P_L(\mathcal{E}_{\mathbb{Z}_{n}}, x)= & \frac{x(x-T)}{(x-m)(x-(T-m))}P_L(K_m, (x-(T-m))) \\ 
   \times  & P_L(\overline{K_{m_1}}\vee \overline{K_{m_2}}, (x-m)),\\
     = & \frac{x(x-T)}{(x-m)(x-(T-m))}(x-(T-m))(x-T)^{m-1} \\
   \times   & P_L(\overline{K_{m_1}}\vee \overline{K_{m_2}}, (x-m)), \\
     = &  x(x-T)^m (x-(T-m_1)^{m_1-1}(x-(T-m_2))^{m_2-1}, \\
\end{split}
\end{equation*} where $m=m_1m_2-1$ and $T=m+m_1+m_2$.\\
Then, 
\begin{equation}
  a(\mathcal{E}_{\mathbb{Z}_{n}})= min\{T-m_1, T-m_2\}= T-max\{m_1, m_2\}.  
\end{equation} Also, by Proposition \ref{alg cnncty of EZn}, 
\begin{equation}
    \kappa(\mathcal{E}_{\mathbb{Z}_{n}})= m+\eta = m+ min\{ m_1,m_2\}.
\end{equation}
\end{enumerate} Thus Equations $(2)$ and $(3)$ guarantees that $a(\mathcal{E}_{\mathbb{Z}_{n}})= \kappa(\mathcal{E}_{\mathbb{Z}_{n}})$.
\end{proof}
\section{Conclusion}

In this paper, we have determined the structure of the subgraph induced by the nonessential ideals of $\mathbb{Z}_{n}$ as the generalized join of certain null graphs. Also, we have proved that the underlying graph in this generalized join is the annihilating ideal graph of $\mathbb{Z}_{n}$ for which $n$ is a product of distinct primes. Further, the various spectra of the induced subgraph of $\mathcal{E}_{\mathbb{Z}_{n}}$ such as adjacency, Laplacian, signless Laplacian, and normalized Laplacian are studied. Also, a characterization of $\mathcal{E}_{\mathbb{Z}_{n}}$ to be  Laplacian integral is obtained in terms of the vertex-weighted Laplacian matrix of annihilating ideal graph of $\mathbb{Z}_{n}$  for $n= \prod_{i=1}^k p_i$. 
Moreover, an upper bound for the spectral radius $b(\mathcal{E}_{\mathbb{Z}_{n}})$  and a characterization for the values of $n$ for which $b(\mathcal{E}_{\mathbb{Z}_{n}})=|V(\mathcal{E}_{\mathbb{Z}_{n}})|$ are estimated. Finally,  the vertex connectivity of $\mathcal{E}_{\mathbb{Z}_{n}}$ is computed and the values of $n$ for which $a(\mathcal{E}_{\mathbb{Z}_{n}})= \kappa(\mathcal{E}_{\mathbb{Z}_{n}})$ are discussed.


\end{document}